\documentclass[12pt]{amsart}

\usepackage{xcolor}
\usepackage[normalem]{ulem} 

\newtheorem{thm}{Theorem}[section]

\newtheorem{lem}[thm]{Lemma}
\newtheorem{prop}[thm]{Proposition}
\theoremstyle{definition}
\newtheorem{defn}[thm]{Definition}
\newtheorem{rem}[thm]{Remark}

\newtheorem{ex}[thm]{Example}
\numberwithin{equation}{section}

\theoremstyle{plain}

\usepackage{amsmath}
\usepackage{amsfonts}
\usepackage{amssymb}
\usepackage{color}
\usepackage{mathtools}
\usepackage{bbm}
\newcommand{\be}{\begin{equation}}
	\newcommand{\en}{\end{equation}}
\newcommand{\Lc}{{\mc L}}
\newcommand{\bei}{\begin{itemize}}
	\newcommand{\eni}{\end{itemize}}
\newcommand{\sip}[2]{\left[{#1}|{#2}\right]}
\newcommand{\rip}[2]{\left({#1}|{#2}\right)}
\newcommand{\ip}[2]{\langle{#1}|{#2}\rangle}

\numberwithin{equation}{section}

\newcommand{\mb}{\mathbb}
\newcommand{\mc}{\mathcal}
\newcommand{\eul}{\mathfrak}
\newcommand{\A}{\eul A}
\newcommand{\Ao}{{\eul A}_{0}}

\newcommand{\D}{\mc D}

\newcommand{\K}{\mc K}

\newcommand{\QAal}{{\mathcal Q}_{\Ao}^\alpha(\A)}

\newcommand{\id}{{\sf e }}

\newcommand{\FF}{\mathfrak F}

\newcommand{\B}{{\eul B}}

\newcommand{\ad}{^{\mbox{\scriptsize $\dag$}}}

\newcommand{\mult}{\,{\scriptstyle \square}\,}

\newcommand{\LDKr}{{\mathcal L}^\ddagger(\mathcal{D},\mathcal{K})}

\newcommand{\LDH}{{\mathcal L}\ad(\D,\H)}
\newcommand{\idop}{{\mb I}}

\def\MM{{\mathfrak M}}

\numberwithin{equation}{section}

\def\H{\mc H}
\def\x{\relax\ifmmode {\mbox{*}}\else*\fi}

\newcommand{\wmult}{\mbox{\raisebox{1pt}{$\scriptscriptstyle{
				\diamond}$}}}
			
\newcommand{\vertiii}[1]{{\left\vert\kern-0.25ex\left\vert\kern-0.25ex\left|#1 
					\right\vert\kern-0.25ex\right\vert\kern-0.25ex\right\vert}}
                    
\definecolor{Magenta}{rgb}{1,0,1}
\definecolor{ao(english)}{rgb}{1,0,0.4}

\newcommand{\BibTeX}{B\kern-0.1emi\kern-0.017emb\kern-0.15em\TeX}
\newcommand{\XYpic}{$\mathrm{X\kern-0.3em\raisebox{-0.18em}{Y}}$-$\mathrm{pic}\,$}

\newcommand{\ed}{\end{document}}
\begin{document}

%
%
%
%
%
%
%
%
%

	\title[]{GNS-like representations of a quasi *-algebra in  Krein spaces}

	\author[G. ~Bellomonte]{Giorgia Bellomonte}

\address{Dipartimento di Matematica e Informatica, Universit\`a degli Studi  di Palermo, Via Archirafi n. 34,  I-90123 Palermo, Italy}
\email{giorgia.bellomonte@unipa.it} 

%
%

\date{July 3, 2026}

\begin{abstract}
In this paper unbounded representations of quasi *-algebras
on an indefinite inner product space are introduced and studied. After introducing the partial *-algebra of unbounded closable
operators sharing a common dense domain in a Krein space,  the
appropriate notion of *-representation in this setting is developed, and
a GNS-type construction for quasi *-algebras based on
$\alpha$-positive sesquilinear forms is established. Furthermore,  the
notion of abstract conditional expectation  is generalized  to quasi *-algebras and
 how {\em admissible} abstract conditional expectations $P$ give rise to representable $P$-forms  is shown. The relation between
$\alpha$-positive forms and $P$-forms is also clarified. Several examples
are provided.
\end{abstract}
\label{page:firstblob}
\maketitle
 \section{Introduction and preliminaries} In the closing years of the twentieth century several works explored
massless and gauge fields by allowing \emph{non-positive} states on the
Borchers field algebra, and by constructing the corresponding quantum
field theories through the Gelfand--Naimark--Segal (GNS) procedure.
The resulting GNS representations naturally led to spaces with
\emph{indefinite metric}, and the mathematical control of these models
required their embedding into Krein spaces. This raises the following
question: is it possible to generalize such a representation to quasi
$*$-algebras by means of sesquilinear forms adapted to an indefinite
metric structure?

To address this question, in Section \ref{sec part *-algebra LDK} the
partial $*$-algebra ${\LDKr}$ of unbounded closable operators sharing a
common dense domain $\D$ in a Krein space $\K$ is first introduced; this provides
the operator framework for the representations considered later. Given
a quasi $*$-algebra $(\A,\A_0)$ with unit $\id$, a class $j(\A)$ of
involutive linear maps $\alpha$ on $\A$, compatible both with the
involution and with the quasi $*$-algebra structure, is then considered.
A subclass of such maps is selected in order to introduce the
corresponding notion of $J_\alpha$-representation. In the same section, the appropriate notion of an $\alpha$-cyclic
$J_\alpha$-\!\! representation in a Krein space is developed. The closure of
such a representation is also introduced and it is proved that it is
again an $\alpha$-cyclic $J_\alpha$-representation of $(\A,\A_0)$ in a
Krein space.

In Section \ref{sec repr in a Krein space}, the class $\QAal$ of
sesquilinear forms $\Phi$ such that the associated form
$$
(a,b)\mapsto \Phi(a,\alpha(b))
$$
is positive is introduced, together with suitable left-invariance and
covariance conditions; see Definition
\ref{defn: alpha left-inv alpha pos sesq form}. This framework includes
genuinely indefinite forms, but does not require $\Phi$ itself to be
indefinite. For forms $\Phi\in\QAal$, under the natural density
assumption on $\lambda^\alpha_\Phi(\A_0)$, where
$\lambda^\alpha_\Phi$ denotes the quotient map associated with
$$
N^\alpha_\Phi:=\{a\in\A:\Phi(a,\alpha(a))=0\},
$$
a GNS-like construction in a Krein space is established.

Finally, in Section \ref{sec abstract expect}, the notion of abstract
conditional expectation is extended to quasi $*$-algebras. Admissible
abstract conditional expectations $P$, characterized by the condition
$$
\alpha_P=2P-I\in j(\A),
$$
are introduced, and the corresponding class of $P$-forms is defined. It
is shown that every $P$-form associated with an admissible abstract
conditional expectation belongs to
$
\mathcal Q^{\alpha_P}_{\A_0}(\A).
$
Conversely, every $\alpha\in j(\A)$ determines the canonical admissible
abstract conditional expectation
$$
P_\alpha=\frac{I+\alpha}{2},
$$
and the forms in $\QAal$ are precisely the $P_\alpha$-forms. Several
examples illustrate the introduced notions.\bigskip\\

In the remainder of this section,  the main definitions and known results that will be needed later are presented. \\

Throughout the paper, for a linear operator $T$ in Hilbert space $\H$, the domain, kernel, range, adjoint, and restriction of $T$ to a subspace $\D$ will be denoted by $D(T)$, $\ker T$, $R(T)$, $T^*$, and $T{\upharpoonright}_{\D}$, respectively.
 Also, $\B(\H)$ will denote the space of bounded linear maps on the Hilbert space $\H$ and $\mathbb{I}$ the identity operator on $\H$.  

 \begin{defn}\cite[Definition 2.1.1]{ait_book}
    A partial *-algebra is a complex vector space $\mathfrak{A}$ with an
involution $x \mapsto x^*$ (i.e. $(x+y)^*=x^*+y^*$, $(\lambda x)^*
=\overline{\lambda}x^*$ ($\lambda\in\mathbb{C}$), $x^{**}=x$) and a subset $\Gamma \subseteq \mathfrak{A}
\times \mathfrak{A}$ such that:

\begin{itemize}
    \item[(i)] $(x,y)\in\Gamma$ iff $(y^*,x^*)\in\Gamma$;
    \item[(ii)] if $(x,y_{1})$ and $(x,y_{2})\in\Gamma$, then $(x,\lambda y_{1}+\mu y_{2})\in\Gamma$, for all $\lambda,\mu\in\mathbb{C}$;
    \item[(iii)] if $(x,y)\in\Gamma$, then there exists an element $x \cdot y \in \mathfrak{A}$
    and  the usual properties of the multiplication hold true:
    $$ x \cdot (y+\lambda z)=x \cdot y + \lambda x \cdot z,\quad (x,y),\, (x,z)\in\Gamma$$ and
$$(x \cdot y)^*=y^* \cdot x^*.$$
An element $\id\in\A$ is called a unit if $\id^*=\id$, $(\id,x)\in\Gamma$, $x \cdot \id=\id  \cdot x=x$ for every $x\in\A$. 
\end{itemize}
\end{defn}

If $(x, y)\in\Gamma$, $x$ is called a left multiplier of $y$ and
$y$ is a right multiplier of $x$ with the notation $x \in L(y)$ and $y \in R(x)$. The universal left- and right-multipliers of $\A$ are indicated, respectively, by $$L\A=\{x\in\A: x\in L(y),\,\, \forall y\in\A\}$$ and 
$$R\A=\{y\in\A: y\in R(x),\,\, \forall x\in\A\}.$$
\begin{ex}
    Let $\D$ be a dense subspace of a Hilbert space $\H$ and let us consider the following family of closable linear operators acting on $\H$ (see \cite[Section 2.1.2]{ait_book}):
\begin{align*}		{\LDH}&=\{X \mbox{ closable}: D(X)=\D;\; D(X^*)\supset \D\} .
\end{align*}
This family is  a partial *-algebra with respect to the following operations: the addition $T+V$, the scalar multiplication $cT$, the involution  defined by		$T^\dag := T^*{\upharpoonright}_{\D}$, the restriction of $T^*$, the adjoint of $T$, to $\D$ and the weak multiplication  $T\mult V:=(T^\dag)^* V$ defined whenever $V\D\subseteq D((T^\dag)^*)$ and $T^\dag\D\subset D(V^*)$.

\end{ex}

A particular case of partial *-algebra is a quasi *-algebra where $\Gamma=\{(a,b)|\,a\in\A_0 \vee b\in  \A_0\}$; its definition is given in full details for the sake of the reader.
\begin{defn}
    A {\em quasi *-algebra} $(\A, \A_0)$ is a pair consisting of a vector space $\A$ and a *-algebra $\A_0$ contained in $\A$ as a subspace and such that
\begin{itemize}
	\item[(i)] $\A$ carries an involution $a\mapsto a^*$ extending the involution of $\A_0$;
	\item[(ii)] $\A$ is  a bimodule over $\A_0$ and the module multiplications extend the multiplication of $\A_0$. In particular, the following associative laws hold:
	\begin{equation}\notag \label{eq_associativity}
		(ca)d = c(ad); \ \ a(cd)= (ac)d, \quad \forall \ a \in \A, \  c,d \in \A_0;
	\end{equation}
	\item[(iii)] $(ac)^*=c^*a^*$, for every $a \in \A$ and $c \in \A_0$.
\end{itemize}

The
\emph{identity} or {\it unit element} of $(\A, \A_0)$, if any, is a necessarily unique element $\id\in \A_0$, such that
$a\id=a=\id a$, for all $a \in \A$.

It will be always supposed that
\begin{align*}
	&ac=0, \; \forall c\in \A_0 \Rightarrow a=0 \\
	&ac=0, \; \forall a\in \A \Rightarrow c=0. 
\end{align*}
Clearly, both these conditions are automatically satisfied if $(\A, \A_0)$ has an identity $\id$.\end{defn}

 \begin{defn} \label{defn_starrepmod} 
Let $(\A,\A_0)$ be a quasi *-algebra with unit $\id$ and $\D_\pi$ a dense subspace
of a certain  Hilbert space $\mathcal{H}$. 
A linear map $\pi$ from $\A$ into  ${\mathcal L}\ad(\D_\pi,\mathcal{H})$ is called  a \emph{*-representation} of  $(\A, \A_0)$,
if the following properties are fulfilled:
\begin{itemize}
\item[(i)]  $\pi(a^*)=\pi(a)^\dagger:=\pi(a)^*{\upharpoonright}_{\D_\pi}, \quad \forall \ a\in \A$;
\item[(ii)] for $a\in \A$ and $c\in \A_0$, $\pi(a)\mult\pi(c)$ is well-defined and \linebreak {$\pi(a)\mult \pi(c)=\pi(ac)$}.
\end{itemize}
For every *-representation $\pi$ of $(\A,\A_0)$,  it is assumed that $\pi(\id)={\idop_{\D_\pi}}$, the
identity operator on  the space $\D_\pi$.\\
 The *-representation $\pi$  is said to be\begin{itemize}
\item[(iii)]  \emph{cyclic} if there  exists $\xi\in\D_\pi$ such that $\pi(\A_0)\xi$ is dense in $\mathcal{H}$ in its  norm topology. In this case $\xi$ is called {\em cyclic vector}.
\end{itemize}
\end{defn}If $\pi$ is a *-representation of $(\A,\A_0)$, then the {\em closure} $\widetilde{\pi}$ of $\pi$ is defined, for each $a\in\A$, as the restriction of the closure $\overline{\pi(a)}$ of the operator $\pi(a)$ to the domain $\widetilde{\D_\pi}$ completion with respect to the graph topology defined by the seminorms $\eta\to \|\eta\|+\|\pi(a)\eta\|$, $a\in\A$, where $\|\cdot\|$ is the norm induced by the inner product  on $\D_\pi$; $\widetilde{\pi}$ is a *-representation of $(\A,\A_0)$, see, e.g.,  \cite[Proposition 7.1.6]{ait_book}.\\ If $\pi=\widetilde{\pi}$, the *-representation $\pi$ is said to be {\em closed}.\\

In \cite[Proposition 2.4.1]{FT_book},  a *-representation of a quasi *-algebra $(\A,\A_0)$ with unit is built in Hilbert space by means of a left *-invariant positive sesquilinear form $\Psi$ which is such that the quotient $\lambda_\Psi(\A_0)=\A_0/N_\Psi$ is dense in the completion $\H_\Psi$ of $\A/N_\Psi$  (with $N_\Psi$ the space of those elements $a\in\A$ such that $\Phi(a,a)=0$) with respect to the norm induced on $\A/N_\Psi$ by the inner product\begin{align*} 	\ip{\lambda_{\Psi}(a)} {\lambda_{\Psi}(b) }_\Psi := \Psi(a, b), \quad \forall \; a,b \in \A.
		\end{align*} More precisely, 
\begin{prop}\cite[Proposition 2.4.1]{FT_book}\label{GNS with a pos l.i. s.f.}
    Let $(\A, \A_0)$ be a quasi *-algebra with unit $\id$ and $\Psi$ a sesquilinear
form on $\A \times\A$ such that  $\Psi$ is positive and left *-invariant (i.e. $\Psi(a, a)\geq 0$, for all $a\in\A$;
and $ \Psi(ax, y) = \Psi(x, a^*y)$, for all $a\in\A$, $x, y \in\A_0$). Then following statements are equivalent:\begin{itemize}
    \item[(i)] $\lambda_\Psi(\A_0)=\A_0/N_\Psi$ is dense in $\H_\Psi$;
\item[(ii)] there exist a Hilbert space $\H_\Psi$, a dense domain $\D_\Psi$ of the Hilbert space $\H_\Psi$ and a closed cyclic *-representation $\pi_\Psi$ in $\mathcal{L}\ad(\D_\Psi,\H_\Psi)$, with cyclic vector $\xi_\Psi$
(in the sense that $\pi_\Psi(\A_0)\xi_\Psi$ is dense in $\H_\Psi$), such that
$$\Psi(a, b) =\ip{\pi_\Psi(a)\xi_\Psi}{\pi_\Psi(b)\xi_\Psi}, \quad \forall a, b \in\A. $$
\end{itemize} 
\end{prop}
\bigskip

Given a Hilbert space  $(\K,\ip{\cdot}{\cdot})$
and a symmetry on it (i.e., an operator $J\in\B(\K)$  such that  $J^2= \idop$ and
$J=J^*$), 
then an indefinite inner product can be defined on  $\K$
by \begin{equation}\label{(a3)}
  \sip{x}{y}   =  \ip{Jx}{y}  ,\ \ 
 x , y \in {\mathcal K}.
 \end{equation} 
A Hilbert space equipped with   an indefinite inner product $\sip{\cdot}{\cdot}$ as in \eqref{(a3)} for some symmetry $J$  
is said to be a {\em Krein space}. 
The Krein space will be often denoted by $\{\K,J\}$ to specify the chosen fundamental symmetry on it.  Denote the positive definite inner product on $\K$ by  \begin{equation}\label{eq def pos inner prod J}
    \rip{\eta}{\zeta}_J=\sip{\eta}{J\zeta}, \quad \forall \eta,\zeta\in\K
\end{equation}
and the induced norm by $$\|\eta\|_J=\sqrt{\rip{\eta}{\eta}_J}=\sqrt{\sip{\eta}{J\eta}}, \quad \eta\in \K.$$ The symbol $\K_
J$ will seldom be used when it is desirable to emphasize that $\K$ is regarded as the Hilbert space endowed with the positive-definite inner product defined in \eqref{eq def pos inner prod J}. Then, \begin{equation*}
    |\sip{\eta}{\zeta}|\leq\|\eta\|_J\|\zeta\|_J,\quad\forall \eta,\zeta\in\K.
\end{equation*} All topological notions in a Krein space 
 will  refer to the norm topology ($\|\cdot\|_J$).

\bigskip

\begin{defn}\cite{AO} Let $\{\K,J\}$ be a  Krein space with indefinite inner product $\sip{\cdot}{\cdot}$ and fundamental symmetry $J$. Let $T:D(T)\to\K$ be a linear operator with $D(T)$ dense in $\K$. Set \begin{multline*}
    D(T^J)=\{\eta\in\K:\, \mbox{there exists }\zeta\in\K \mbox{ such that }\\\sip{T\xi}{\eta}=\sip{\xi}{\zeta},\, \forall\xi\in D(T)\}.
\end{multline*} The well-defined linear operator $$T^J\eta:=\zeta, \quad\eta\in D(T^J)$$ is said to be the {\em $J$-adjoint $T^J$ of $T$}.
\end{defn}

Given a linear operator $T$ in Hilbert space $\K_J$,  $T^*$ will indicate as usual its (Hilbertian) adjoint.\\
    
Now some known facts on densely defined linear operators in Krein spaces that will be used throughout the paper are gathered (see either \cite[Sect. C]{Jakobczyk86} or  \cite[Ch. VI]{BOGNAR}).

\begin{lem}\label{lem 3.1}
 Let  $\{\K,J\}$ be  a  Krein space  with fundamental symmetry $J$. Let $\D\subset\K$ be dense in $\K$ and a linear operator $T:\D\to \K$ and  $T^J$ its $J$-adjoint. Then,\begin{itemize}
     \item[$i$)] $T^J=JT^*J$,
        \item[$ii)$]  $T^J$ is closed.\end{itemize}Moreover,  \begin{itemize}\item[$iii)$] $T$ is closable iff $D(T^J)$ is dense in $\K$;
     \item[$iv)$]  if  $D(S)\subset\K$ is dense and $S\subset T$, then $T^J\subset S^J$; 
     \item[$v)$] if $T$ is closable then $T^{JJ}=(T^J)^J$ exists and in this case  $T^{JJ}=\overline{T}$. 
         \end{itemize}
\end{lem}
\begin{proof}  
    $i)$   By \cite[Ch. VI, Lemma 2.1]{BOGNAR}. $ii)$ By \cite[Ch. VI, Theorem 2.2]{BOGNAR}.
 $iii)$ By \cite[Proposition 2.1]{Jakobczyk86}.
 $iv)$ By \cite[Ch. VI, Theorem 2.6]{BOGNAR}.
 $v)$ By \cite[Ch. VI, Theorem 2.6]{BOGNAR} combined with $iv)$.
 \end{proof}

The following characterization, which is essentially a reformulation of standard results on fundamental symmetries in Krein spaces (see \cite[Chapter I]{BOGNAR}) will be used later.

\begin{thm}\label{thm extendability of the inner prod} Let $(\D,\sip{\cdot}{\cdot})$ be a non degenerate indefinite inner product space.
Suppose that there exists an operator
$J:\D\to\D$ satisfying:
\begin{itemize}
    \item[(1)] $J^2 = \idop_\D$;
    \item[(2)] $\sip{J\xi}{\eta} = \sip{\xi}{J\eta}$ for all $\xi,\eta\in\D$;
    \item[(3)] $(\xi,\eta)_J := \sip{\xi}{J
    \eta}$ is a positive definite inner product on $\D$.
\end{itemize}
Let $\K=\widetilde{\D}^{\|\cdot\|_J}$ be the Hilbert space completion of $(\D,(\cdot,\cdot)_J)$ with respect to the norm $\|\cdot\|_J=\sqrt{(\cdot,\cdot)_J}$. Then the following statements are equivalent:
\begin{enumerate}
    \item[(i)] there exists an indefinite inner product $\sip{\cdot}{\cdot}'$ on $\K$
    extending $\sip{\cdot}{\cdot}$  such that $(\K,\sip{\cdot}{\cdot}')$ is a Krein space;
    \item[(ii)] the operator $J$ is bounded with respect to the norm
    $\|\cdot\|_J$, i.e.,
    $$\exists\,C>0:\quad \|J\eta\|_J \le C\|\eta\|_J,\qquad \forall\,\eta\in\D.
    $$
\end{enumerate}
If one of the previous statements is true, $J$ extends uniquely to a bounded self-adjoint involution
$\overline{J}:\K\to\K$, and the extended indefinite inner product is given by
$$ \sip{\xi}{\eta}':= (\xi,\overline{J}\eta)_J,\qquad \xi,\eta\in\K.$$
Moreover, $(\K,\sip{\cdot}{\cdot}')$ is a Krein space with fundamental symmetry
$\overline{J}$.
\end{thm}

\section{(Closed) $\alpha$-cyclic $J_\alpha$-representations in Krein spaces}\label{sec part *-algebra LDK}

In order to define a suitable notion of *-representation of quasi *-algebras in Krein spaces, the operator framework needed to it has to be first introduced.

\begin{defn}\label{defn LDKr}
    Let $\D$ be a dense subspace of the Krein space $\{\K, J\}$  with indefinite inner product $\sip{\cdot}{\cdot}$ and fundamental symmetry $J$. Let us consider the following families of linear operators acting on $\D$:
\begin{equation*}
    {\LDKr}=\{T \mbox{ closable linear operator in } \K: D(T)=\D;\, D(T^J)\supset \D\}
\end{equation*}
    \begin{equation*}
       {\Lc^\ddagger(\D)}=\{T\in \LDKr: T\D\subset \D; \; T^J\D\subset \D\}. 
    \end{equation*}
\end{defn}

 The set  ${\LDKr}$ is a vector space with respect to the  addition $T+V$ and the scalar multiplication $cT$,  if $T,V\in{\LDKr} $ and $c\in\mathbb{C}$. \\

 The following proposition is a direct consequence of Lemma \ref{lem 3.1}.

\begin{prop}\label{prop LDK}Let $\D$ be a dense subspace of the Krein space $\{\K, J\}$  with indefinite inner product $\sip{\cdot}{\cdot}$ and fundamental symmetry $J$.
      The vector space ${\LDKr}$ is a partial *-algebra with respect to the following operations: the addition $T+V$, the scalar multiplication $cT$, the involution $\ddagger:T\to T^\ddag=T^J{\upharpoonright}_{\D}$ and the weak multiplication  $T\wmult V:=(T^\ddag)^J V$ defined whenever $V\D\subset D((T^\ddag)^J)$ and $T^\ddag\D\subset D(V^J)$ and has the identity operator  $\idop{\upharpoonright}_\D
      $ as unit.     The set $\Lc^\ddagger(\D)$ is a  *-algebra with respect to the usual multiplication $XY$, for all $X,Y\in\Lc^\ddagger(\D)$ (which here coincides with the weak multiplication $X\wmult Y$) and with the involution $\ddagger$.
\end{prop}

  Let $\alpha:\A\to\A$ be a map. Throughout this paper, the set of fixed points of $\alpha$ will be denoted by
$$
F_\alpha(\A)=\{a\in\A:\alpha(a)=a\},
$$
and, for every subset $\A_0\subset\A$, let
$$
F_\alpha(\A_0)=F_\alpha(\A)\cap\A_0.
$$

\begin{defn}\label{defn: alpha map}
    Let $(\A,\A_0)$ be a quasi *-algebra with unit $\id$.   The set of all those linear maps $\alpha:\A\to\A$ such that
    \begin{itemize}
		    \item[($\alpha$.1)] $\alpha(\id)=\id$;
            \item[($\alpha$.2)] $\alpha(a)^*=\alpha(a^*)$, for all $a\in\A$;
            \item[($\alpha$.3)] $\alpha(\A_0)\subseteq\A_0$, $\alpha(ax)=\alpha(a)\alpha(x)$, for all $a\in\A$, $x\in\A_0$; 
            \item[($\alpha$.4)] $\alpha$ is a linear involution on $\A$: $\alpha^2(a)=a$, for all $a\in\A$.
		\end{itemize} will be denoted by $j(\A)$.
                \end{defn} 
                
                \begin{rem} \label{rem fixed point of alpha in JA}If $\alpha\in j(\A)$, then the set $F_\alpha(\A)$ 
is a non empty *-invariant linear subspace of $\A$ (in particular, $\id\in F_\alpha(\A_0)\subset F_\alpha(\A)$) and 
the set $F_\alpha(\A_0)$ is  a unital *-subalgebra of $\A_0$.
\end{rem}

\begin{ex}\label{ex alpha}
    Let us consider  $(\A,\A_0)=(L^2([0,1]),L^\infty([0,1]))$ where multiplication is the pointwise multiplication whenever one of the
	two factors belongs to $\A_0$, and the involution is given by 	$
	f^*=\overline f.
	$
	Then $(L^2([0,1]),L^\infty([0,1]))$ is a quasi *-algebra with unit $\id=1$ the constant function on $[0,1]$. Consider the reflection $\alpha:L^2([0,1])\to L^2([0,1])$ $$\alpha(f)(t)=Rf(t)=f(1-t).$$ It is easy to see that $\alpha\in j(L^2([0,1]))$. 
\end{ex}\medskip

    \begin{defn} \label{defn-J-repre} 
Let $(\A,\A_0)$ be a quasi *-algebra with unit $\id$ and $(\D,\sip{\cdot}{\cdot})$ be a non degenerate indefinite inner product space.
The set  of those $\alpha\in j(\A)$ such that there exists  an operator $J_\alpha:\D\to\D$      such that \begin{itemize}
    \item[$(j1)$] $(J_\alpha)^2=\idop_{\D}$;
    \item[$(j2)$]  $\sip{J_\alpha\zeta}{\eta}=\sip{\zeta}{J_\alpha\eta}$ for all $\zeta,\eta\in\D$, \item[$(j3)$]  $\rip{\eta}{\zeta}_{J_\alpha}:=\sip{\eta}{J_\alpha\zeta}$, $\eta,\zeta\in\D$, is a positive definite inner product with associated norm $\|\cdot\|_{J_\alpha}$,
    \end{itemize} will be  denoted by   $j_\D(\A)$. \end{defn}

\begin{rem}\label{rem boundedness of J(alpha)}
    If $\alpha\in j_{\D}(\A)$, then the operator $J_\alpha:\D\to\D$ in Definition  \ref{defn-J-repre} is bounded on $\D$ (more precisely it is isometric):
        \begin{align*}        \|J_\alpha\eta\|_{J_\alpha}^2&=\rip{J_\alpha\eta}{J_\alpha\eta}_{J_\alpha}:=\sip{J_\alpha\eta}{J_\alpha^2\eta}\\&=\sip{J_\alpha\eta}{\eta}=\sip{\eta}{J_\alpha\eta}=\|\eta\|_{J_\alpha}^2,
    \end{align*}  therefore 
        by Theorem \ref{thm extendability of the inner prod}, there exists a (non degenerate) inner product $\sip{\cdot}{\cdot}'$ on $\K$ which is the extension  of $\sip{\cdot}{\cdot}$  to $\K$ .
\end{rem}
\medskip

   The partial *-algebra structure described in Proposition
\ref{prop LDK} allows us to formulate the following notion of representation.
    \begin{defn}\label{defn J alpha repres} Let $(\A,\A_0)$ be a quasi *-algebra with unit $\id$, $(\D_\pi,\sip{\cdot}{\cdot})$ be a non degenerate indefinite inner product space and $\alpha\in j_{\D_\pi}(\A)$. Let  $\K_\pi$ be the 
    completion of $\D_\pi$ with respect to the norm $\|\cdot\|_{J_\alpha}$ 
    defined as in $(j3)$ and let $\sip{\cdot}{\cdot}_\pi$ be
 the extension of $\sip{\cdot}{\cdot}$  to $\K_\pi$ (therefore $\{\K_\pi,\sip{\cdot}{\cdot}_\pi\}$ is a Krein space).         A linear map $\pi$ from $\A$ into ${\mathcal L}^\ddagger(\D_\pi,\K_\pi)$  is called a {\em $J_\alpha$-representation of $(\A,\A_0)$} if:
\begin{itemize}
\item[(r1)]  for every $a\in\A$,  $\pi(a^*)=\pi(\alpha(a))^\ddagger:=\pi(\alpha(a))^{J_\alpha}{\upharpoonright}_{\D_\pi}$;
\item[(r2)] for every $a\in \A$ and $c\in \A_0$, $\pi(a)\wmult\pi(c)$ is well-defined  and  $\pi(a)\wmult \pi(c)=\pi(ac)$;
    \item[(r3)] for every $a\in\A$, $\eta\in\D_\pi$, $\pi(\alpha(a))\eta=J_\alpha\pi(a)J_\alpha\eta$.
\end{itemize}
It is assumed that for every $J_\alpha$-representation $\pi$ of $(\A,\A_0)$,  $\pi(\id)=\idop_{\D_\pi}$, the
identity operator on  the space $\D_\pi$.\\
The $J_\alpha$-representation of $(\A,\A_0)$ $\pi$ is called  {\em $\alpha$-cyclic}  if it is a $J_\alpha$-representation such that 
\begin{itemize}
    \item[(r4)]   there exists $\xi_\pi\in\D_\pi$ such that $\pi(\A_0)\xi_\pi$ is dense in $\K_\pi$ in its norm topology $\|\cdot\|_{J_\alpha}$ and $J_\alpha\xi_\pi=\xi_\pi$. The vector $\xi_\pi$ is called {\em $\alpha$-cyclic}.
\end{itemize}\end{defn}

Since the aim is to obtain a closed $J_\alpha$-representation of $(\A,\A_0)$ in Krein spaces, the closure $\widetilde{\pi}$ of a $J_\alpha$-representation $\pi$ in a Krein space is introduced, analogously to the Hilbert space case, and it is proved that $\widetilde{\pi}$ is still a $J_\alpha$-representation.

Let $\alpha\in j_{\D_\pi}(\A)$, $\pi$ be a $J_\alpha$-representation of $(\A,\A_0)$ and  $t_\pi$ be the graph topology induced by the family of seminorms \begin{equation}
    \label{eq: seminorms}
p_a(\eta)=\|\eta\|_{J_\alpha}+\|\pi(a)\eta\|_{J_\alpha},\quad a\in\A\end{equation} where $\|\cdot\|_{J_\alpha}$ is the norm  induced by the definite positive inner product defined in $(j3)$ on $\D_\pi$. Put $$\widetilde{\D_\pi}:=\widetilde{\D_\pi}^{t_\pi}$$ i.e. the completion of $\D_\pi$ with respect to $t_\pi$. The canonical embedding of $\widetilde{\D_\pi}$ in $\K_\pi$ is well-defined and injective.  \\Since every operator $\pi(a)$, $a\in\A$, is closable, its closure $\overline{\pi(a)}$ can be considered. Then $\widetilde{\D_\pi}$ can be identified with a subset of $\bigcap_{a\in\A}D(\,\overline{\pi(a)}\,)$. This is because for every $\eta\in \widetilde{\D_\pi}$ there exists a net $(\eta_r)\subset\D_\pi$ such that  both $(\eta_r)$ and $(\pi(a)\eta_r)$ for every $a\in\A$, are Cauchy in the norm topology.
For every $a\in\A$, since $\pi(a)$ is closable, the limit of $(\pi(a)\eta_r)$ does not depend on the chosen net $(\eta_r)$  whose limit is $\eta$; moreover, the limit of  the couple $(\eta_r,\pi(a)\eta_r)$
 belongs to the graph of $\overline{\pi(a)}$.

Therefore,
$$
\widetilde{\D_\pi}\subseteq\bigcap_{a\in\A}D(\,\overline{\pi(a)}\,),
$$
and one can define $\widetilde{\pi}(a)\eta=\lim_r\pi(a)\eta_r$; in this way, similarly to the Hilbert space case, the {\em closure} $\widetilde{\pi}$ of $\pi$ is defined, for each $a\in\A$, by
\begin{equation}
     \label{eq closed J alpha repr}
     \widetilde{\pi}(a):=\overline{\pi(a)}{\upharpoonright}_{\widetilde{\D_\pi}},
\end{equation}
that is, as the restriction of the closure $\overline{\pi(a)}$ of the operator $\pi(a)$ to the domain $\widetilde{\D_\pi}$.

 \begin{lem}\label{lem pi tilde in L(D tilde, K)}Let $(\A,\A_0)$ be a quasi *-algebra with unit $\id$, $(\D_\pi,\sip{\cdot}{\cdot})$ be a non degenerate indefinite inner product space and $\alpha\in j_{\D_\pi}(\A)$. 
 If $\pi$ is a $J_\alpha$-representation of $(\A,\A_0)$, then
      $\widetilde{\pi}$ is linear and maps $\A$ into ${\mathcal L}^\ddagger(\widetilde{\D_\pi},\K_\pi)$ and $\A_0$ into ${\mathcal L}^\ddagger(\widetilde{\D_\pi})$.
 \end{lem}\begin{proof}  Let $\pi$ be a $J_\alpha$-representation of $(\A,\A_0)$ into ${\mathcal L}^\ddagger(\D_\pi,\K_\pi)$ and consider the map defined as in \eqref{eq closed J alpha repr}.

It remains to prove that $\widetilde{\pi}(a)\in{\mathcal L}^\ddagger(\widetilde{\D_\pi},\mathcal{K}_\pi)$ for every $a\in\A$ and that $\widetilde{\pi}(x)\in{\mathcal L}^\ddagger(\widetilde{\D_\pi})$ for every $x\in\A_0$. By definition, $\widetilde{\pi}$ is clearly linear. Moreover, all the operators $\widetilde{\pi}(a)$ have the common domain
$
D(\widetilde{\pi}(a))=\widetilde{\D_\pi}
$
and are closable, since $\widetilde{\pi}(a)$ is a restriction of $\overline{\pi(a)}$ for every $a\in\A$.

Now, notice that $
J_\alpha\widetilde{\D_\pi}\subseteq\widetilde{\D_\pi}.
$
Indeed, let $\eta\in\widetilde{\D_\pi}$. Then there exists a $t_\pi$-Cauchy net $(\eta_r)\subset\D_\pi$ representing $\eta$. Since $J_\alpha$ is continuous, the net $(J_\alpha\eta_r)$ is norm-Cauchy. Furthermore, for every $a\in\A$, one has
$
\pi(a)J_\alpha\eta_r=J_\alpha\pi(\alpha(a))\eta_r$, for all $r$.
As $(\pi(\alpha(a))\eta_r)$ is norm-Cauchy and $J_\alpha$ is continuous, it follows that $(\pi(a)J_\alpha\eta_r)$ is also norm-Cauchy. Hence, the net $(J_\alpha\eta_r)\subset\D_\pi$ is $t_\pi$-Cauchy. Consequently,
$
J_\alpha\eta\in\widetilde{\D_\pi},
$
and therefore
$
J_\alpha\widetilde{\D_\pi}\subseteq\widetilde{\D_\pi}$.

Besides, it is $\widetilde{\D_\pi}\subseteq D(\widetilde{\pi}(a)^{J_\alpha})$ for every $a\in\A$, in fact once fixed $a\in\A$ it is enough to show that for every  $\zeta\in \widetilde{\D_\pi}$ there exists a vector
$\xi\in \K_\pi$ such that
$$
\sip{\widetilde{\pi}(a)\eta}{\zeta}
=
\sip{\eta}{\xi},
\qquad \forall \eta\in \widetilde{\D_\pi}.
$$
Let $(\zeta_s)\subset \D_\pi$ be a
$t_\pi$-Cauchy net representing $\zeta\in\widetilde{\D_\pi}$, fix  $\eta\in\widetilde{\D_\pi}$ and consider a $t_\pi$-Cauchy net  $(\eta_r)$ representing it.
Since $\pi$ is a $J_\alpha$-representation, by $(r1)$  it is
$$
\sip{\pi(a)\eta_r}{\zeta_s}
=
\sip{\eta_r}{ \pi(\alpha(a^*))\zeta_s},\quad \forall r,s.
$$

Passing  on both sides to the limits first in $r$ and then in $s$ and using the continuity of the indefinite inner
product with respect to the Hilbert norm $\|\cdot\|_{J_\alpha}$, one obtains
$$
\sip{\widetilde{\pi}(a)\eta}{ \zeta}
=
\sip{\eta}{\widetilde{\pi}(\alpha(a^*))\zeta}.
$$
This is true for every $\eta\in\widetilde{\D_\pi}$, therefore $
\zeta\in D(\widetilde{\pi}(a)^{J_\alpha})
$ 
and by definition of $J_\alpha$-adjoint it is 
\begin{equation}\label{eq widetilde J adjoint}
\widetilde{\pi}(a)^{J_\alpha}\zeta
=
\widetilde{\pi}(\alpha(a^*))\zeta.    
\end{equation} Hence $
\widetilde{\D_\pi}
\subseteq
D(\widetilde{\pi}(a)^{J_\alpha})$, for all $a\in \A$.

   Now, $
\widetilde{\pi}(\mathcal \A_0)\widetilde{\D_\pi}
\subseteq
\widetilde{\D_\pi}
$, i.e.,  $\widetilde{\pi}(x)\widetilde{\D_\pi}\subseteq\widetilde{\D_\pi}$ for all $x\in\A_0$. Indeed, let $\eta\in \widetilde{\D_\pi}$, then there exists a $t_\pi$-Cauchy net $(\eta_r)\subset \D_\pi$
representing $\eta$. This means that $(\eta_r)$ and
$(\pi(a)\eta_r)$ are Cauchy in the norm $\|\cdot\|_{J_\alpha}$,
for every $a\in\mathcal \A$.

Fix $x\in\mathcal \A_0$ and put
$$
\zeta_r:=\pi(x)\eta_r\in \D_\pi .$$

The net $(\zeta_r)$ is $t_\pi$-Cauchy. Indeed, $(\zeta_r)=(\pi(x)\eta_r)$ is  norm-Cauchy, because
$x\in\mathcal \A_0\subseteq\A$. Moreover, for every
$a\in\mathcal \A$,
$$
\pi(a)\zeta_r
=
\pi(a)\wmult\pi(x)\eta_r
=
\pi(ax)\eta_r,
$$
since $x\in\mathcal \A_0$. Since $ax\in\mathcal \A$ and
$(\eta_r)$ is $t_\pi$-Cauchy, also the net
$(\pi(ax)\eta_r)$ is Cauchy in $\|\cdot\|_{J_\alpha}$.

Therefore $(\zeta_r)$ defines an element of
$\widetilde{\D_\pi}$. Under the canonical embedding of
$\widetilde{\D_\pi}$ into $\K_\pi$, this element is
$$
\lim_r \pi(x)\eta_r=\widetilde{\pi}(x)\eta .
$$
Hence $\widetilde{\pi}(x)\eta\in \widetilde{\D_\pi},
$
and therefore
$
\widetilde{\pi}(\mathcal \A_0)\widetilde{\D_\pi}
\subseteq
\widetilde{\D_\pi}.
$

Finally, let $x\in\A_0$ and $\eta\in\widetilde{\D_\pi}$. By
\eqref{eq widetilde J adjoint},
$$
\widetilde{\pi}(x)^{J_\alpha}\eta
=
\widetilde{\pi}(\alpha(x^*))\eta.
$$
Since $\alpha(x^*)\in\A_0$, 
$$
\widetilde{\pi}(\alpha(x^*))\eta\in\widetilde{\D_\pi}.
$$
Hence $
\widetilde{\pi}(x)^{J_\alpha}\eta\in\widetilde{\D_\pi},
$
and therefore
$
\widetilde{\pi}(x)^{J_\alpha}\widetilde{\D_\pi}
\subseteq
\widetilde{\D_\pi}.
$
 \end{proof}

\begin{prop}\label{prop closed repr} Let $(\A,\A_0)$ be a quasi *-algebra with unit $\id$, $(\D_\pi,\sip{\cdot}{\cdot})$ be a non degenerate indefinite inner product space and $\alpha\in j_{\D_\pi}(\A)$. 

 If $\pi$ is a $J_\alpha$-representation of $(\A,\A_0)$, then  
    $\widetilde{\pi}$ is a closed $J_\alpha$-representation. If $\xi_\pi$ is $\alpha$-cyclic for $\pi$ then it is $\alpha$-cyclic for $\widetilde{\pi}$.
\end{prop}
 \begin{proof} 
 Let $\pi$ be a $J_\alpha$-representation of $(\A,\A_0)$ into ${\mathcal L}^\ddagger(\D_\pi,\K_\pi)$ and consider its closure defined as in \eqref{eq closed J alpha repr} on $\widetilde{\D_\pi}$. 
 
     In order to prove that $\widetilde{\pi}$ is a $J_\alpha$-representation, by Lemma \ref{lem pi tilde in L(D tilde, K)}, it remains to check properties $(r1)-(r3)$ in Definition \ref{defn J alpha repres}. Let $\eta\in \widetilde{\D_\pi}$, then there exists a $t_\pi$-Cauchy net $(\eta_r)\subset \D_\pi$
 representing $\eta$. \\ To prove that $(r3)$ holds for $\widetilde{\pi}$, it suffices to observe that \begin{equation}
    \label{eq on D pi}
\pi(\alpha(a))\eta_r=J_\alpha\pi(a)J_\alpha\eta_r,\quad \forall r\end{equation} and by the continuity of $J_\alpha$ and since $
J_\alpha\widetilde{\D_\pi}
\subseteq
\widetilde{\D_\pi}
$, one can pass to the limits on both sides of \eqref{eq on D pi} to get $\widetilde{\pi}(\alpha(a))\eta=J_\alpha\pi(a)J_\alpha\eta$.\\

         To prove $(r2)$, it has to be first shown that
$\widetilde{\pi}(a)\wmult\widetilde{\pi}(x)$ is well-defined, for
$a\in \A$ and $x\in \A_0$. This amounts to proving that
$$
\widetilde{\pi}(x)\widetilde \D_\pi
\subset
D((\widetilde{\pi}(a)^\ddagger)^{J_\alpha})
$$
and
$$
\widetilde{\pi}(a)^\ddagger\widetilde \D_\pi
\subset
D(\widetilde{\pi}(x)^{J_\alpha}).
$$
The first inclusion follows from Lemma \ref{lem pi tilde in L(D tilde, K)}, since
$$
\widetilde{\pi}(x)\widetilde \D_\pi\subset \widetilde \D_\pi
\subset D((\widetilde{\pi}(a)^\ddagger)^{J_\alpha}).
$$

Let now $\zeta\in \widetilde \D_\pi$. Since by \eqref{eq widetilde J adjoint}
\begin{equation*}
    \widetilde{\pi}(b)^\ddagger=\widetilde{\pi}(b)^{J_\alpha}{\upharpoonright}_{\widetilde{\D_\pi}}=\widetilde{\pi}(\alpha(b^*)), \quad \forall b\in\A,
\end{equation*} 
one has to prove that
$$
\widetilde{\pi}(\alpha(a^*))\zeta
\in D(\widetilde{\pi}(x)^{J_\alpha}).
$$
For every $\eta\in \widetilde \D_\pi$, using Lemma \ref{lem pi tilde in L(D tilde, K)} and the fact
that $x\in \A_0$, this yields
$$
\begin{aligned}
\sip{\widetilde{\pi}(x)\eta}{\widetilde{\pi}(\alpha(a^*))\zeta}
&=
\sip{\widetilde{\pi}(ax)\eta}{ \zeta} \\
&=
\sip{\eta}{ \widetilde{\pi}(\alpha((ax)^*))\zeta}.
\end{aligned}
$$
Hence, by the definition of the $J_\alpha$-adjoint,
$$
\widetilde{\pi}(\alpha(a^*))\zeta
\in D(\widetilde{\pi}(x)^{J_\alpha}).
$$
Thus the weak multiplication
$\widetilde{\pi}(a)\wmult\widetilde{\pi}(x)$ is well-defined.

It remains to prove that
$$
\widetilde{\pi}(a)\wmult\widetilde{\pi}(x)=\widetilde{\pi}(ax).
$$
Let $\eta\in \widetilde \D_\pi$ and let
$(\eta_r)\subset \D_\pi$ be a $t_\pi$-Cauchy net representing
$\eta$. Since $\pi$ is a $J_\alpha$-representation,
$$
\pi(a)\wmult\pi(x)\eta_r=\pi(ax)\eta_r, \quad \forall r.
$$
Passing to the limit, it is obtained
$$
\widetilde{\pi}(a)\wmult\widetilde{\pi}(x)\eta
=
\widetilde{\pi}(ax)\eta.
$$
Therefore
$$
\widetilde{\pi}(a)\wmult\widetilde{\pi}(x)=\widetilde{\pi}(ax).
$$
 To prove $(r1)$, observe that, by Lemma \ref{lem pi tilde in L(D tilde, K)}, the restriction $\widetilde{\pi}(a)^\ddagger$ of $\widetilde{\pi}(a)^{J_\alpha}$ to $\widetilde{\D_\pi}$ coincides with $\widetilde{\pi}(\alpha(a^*))$. Then it is also $$\widetilde{\pi}(\alpha(a))^{J_\alpha}{\upharpoonright}_{\widetilde{\D_\pi}}=\widetilde{\pi}(\alpha^2(a^*))=\widetilde{\pi}(a^*)$$ hence $\widetilde{\pi}(\alpha(a))^\ddagger=\widetilde{\pi}(\alpha^2(a^*))=\widetilde{\pi}(a^*)$. Hence $(r1)$ is proved. Moreover  $\widetilde\pi$ is closed by construction.\\  Finally, if $\xi_\pi\in\D_\pi$ is $\alpha$-cyclic for $\pi$, then $\pi(\A_0)\xi_\pi$ is dense in $\K_\pi$ in its norm topology $\|\cdot\|_{J_\alpha}$ and $J_\alpha\xi_\pi=\xi_\pi$. It is clearly  $\xi_\pi\in\D_\pi\subseteq \widetilde{\D_\pi}$, and  since $\pi(\A_0)\xi_\pi\subseteq\widetilde{\pi}(\A_0)\xi_\pi$, it follows that $\widetilde{\pi}(\A_0)\xi_\pi$ is dense  in $\K_\pi$ in its norm topology $\|\cdot\|_{J_\alpha}$ too, therefore $\xi_\pi$ is $\alpha$-cyclic also for $\widetilde{\pi}$.
 \end{proof}

\section{GNS-like construction  for quasi *-algebras in a Krein space via  sesquilinear forms}\label{sec repr in a Krein space}

 In this section a class of sesquilinear forms that provide the basis for the Krein-space GNS construction is introduced. The key point is not to
require $\Phi$ itself to be positive, but rather to assume positivity
of the  form
$$
(a,b)\mapsto \Phi(a,\alpha(b)).
$$
This allows one to construct a Hilbert completion together with an
indefinite inner product whose fundamental symmetry is induced by
$\alpha\in j(\A)$.
\begin{defn}\label{defn: alpha left-inv alpha pos sesq form}
 Let $(\A,\A_0)$ be a quasi *-algebra with unit $\id$, $\alpha\in j(\A)$. The set of all  forms  $\Phi:\A\times \A\to\mathbb{C}$  that are sesquilinear, i.e., \begin{align*}
\Phi(\mu_1a+\mu_2b,\nu_1c+\nu_2d)&=\mu_1\overline{\nu_1}\Phi(a,c)+\mu_1\overline{\nu_2}\Phi(a,d)\\&+\mu_2\overline{\nu_1}\Phi(b,c)+\mu_2\overline{\nu_2}\Phi(b,d), 
       \end{align*} $ a,b,c,d\in\A$, $\mu_1,\mu_2,\nu_1,\nu_2\in\mathbb{C}$ and such that\begin{itemize}\item[q.1)] {\em $\alpha$-positive}, i.e., $\Phi(a, \alpha(a))\geq0$, for all $a\in\A$;
                   \item[q.2)]  {\em $\alpha$-left-invariant}, i.e., $\Phi(ax,\alpha(y))=\Phi(x,\alpha(a^*y)),\,\,\forall a\in\A$, $x,y\in\A_0$;     
     \item[q.3)] {\em covariant} $\Phi(a,b)=\Phi(\alpha(a),\alpha(b))$, for all $a,b\in\A$.
    \end{itemize} will be denoted by $\QAal$.\end{defn}

\begin{ex}\label{ex form in Qaa}
    Let us consider both the quasi *-algebra and the $\alpha\in j(\A)$ of the Example \ref{ex alpha}. Consider the form $\Phi:L^2([0,1])\times L^2([0,1])\to \mathbb{C}$ defined as:$$\Phi(f,g)=\int_0^1\overline{f(t)}g(1-t)dt,\quad f,g\in L^2([0,1]);$$ $\Phi$ is well-defined since if $f,g\in L^2([0,1])$, then $\overline{f(t)}g(1-t)=\overline{f}\alpha(g)(t)\in L^1([0,1])$. It is easy to prove that $\Phi$ is an {\em indefinite} sesquilinear form in ${\mathcal Q}_{L^\infty([0,1])}^\alpha(L^2([0,1]))$.
   \end{ex}

\begin{lem}\label{lem: CSI} Let $(\A,\A_0)$ be a quasi *-algebra with unit $\id$, $\alpha\in j(\A)$ and  $\Phi$ be an  $\alpha$-positive sesquilinear form such that $q.3)$ holds, then the following Cauchy-Schwarz inequality holds
\begin{equation}\label{eq: CS-in H_Phialpha}
|\Phi(a,b)|\leq\Phi( a,\alpha(a))^{1/2}\Phi(b, \alpha(b))^{1/2}, \quad \forall  a,b \in \A.
\end{equation}\end{lem}\begin{proof} Since $\Phi(\cdot,\alpha(\cdot))$ is semidefinite positive, by both $q.3)$ and $(\alpha.4)$ it follows that 
            \begin{align*}
|\Phi(a,b)|&=|\Phi(\alpha(a),\alpha(b))|\\&\leq\Phi(\alpha(a), \alpha(\alpha(a)))^{1/2}\Phi(b, \alpha(b))^{1/2}\\&=\Phi(\alpha(a), a)^{1/2}\Phi(b, \alpha(b))^{1/2}, \quad \forall  a,b \in \A.
        \end{align*} The thesis follows from the observation  that, by the $\alpha$-positivity of $\Phi$, the positive form $\Phi(\cdot,\alpha(\cdot))$ is hermitian.        
    \end{proof}

   In particular \eqref{eq: CS-in H_Phialpha} holds for $\Phi\in\QAal$, with $\alpha\in j(\A)$.\\

\begin{rem}     \label{rem relationship bet adjoints} 
   Let $\alpha\in j(\A)$. 
In the following $T^*_\alpha$ will denote the (Hilbertian) adjoint of $T$ with respect to the $\alpha$-inner product $\Phi_\alpha$ (see $(j3)$). It follows that $T^J = J T^*_\alpha J$.
\end{rem}
\bigskip
Now consider $\alpha\in j(\A)$, $\Phi\in\QAal$  and put
		$$ N_\Phi^\alpha := \big\{ a \in \A : \Phi(a,\alpha(a))= 0\big\};$$
        by Lemma \ref{lem: CSI} $$ N_\Phi^\alpha := \big\{ a \in \A : \Phi(a,\alpha(b))=0, \ \forall 		b \in \A \big\},
		$$
		so $ N_\Phi^\alpha$ is a subspace of $\A$ with the property that, if  
			\begin{equation}\label{eq: property of N}
			  a\in \A,\,\, x \in N_\Phi^\alpha \cap \A_0  \ \text{ then } \Phi(ax,\alpha(y))=0\quad \forall y\in\A_0
			\end{equation} in fact: \begin{align}\label{eq for the well-definiteness}
\Phi(ax,\alpha(y))&=\Phi(x,\alpha(a^*y))=0, \quad \forall \ y \in \A_0.\end{align} 
		As usual  the quotient map $\lambda_{\Phi}^\alpha : \A \rightarrow \A/N_\Phi^\alpha$ is considered. For each $ a
		\in \A$,  $\lambda_{\Phi}^\alpha(a)$ will denote  the corresponding coset of $\lambda_\Phi^\alpha(\A)=
		\A/N_\Phi^\alpha$, which contains $a$.  The space  $\lambda_\Phi^\alpha(\A)$ is a pre-Hilbert space with respect to a definite positive inner product:    
    \begin{equation}
        \label{eq def pos inner prod on lambda A}\rip{\lambda_\Phi^\alpha(a)}{\lambda_\Phi^\alpha(b)}_\alpha=\Phi(a,\alpha(b)),\quad a,b\in\A
    \end{equation}with associated norm $\|\cdot\|_\alpha$. On $\lambda_\Phi^\alpha(\A)$  a non degenerate indefinite inner product  
    can be defined:\begin{equation}
        \label{eq indef inner prod on lambda A}\sip{\lambda_\Phi^\alpha(a)}{\lambda_\Phi^\alpha(b)}=\Phi(a,b),\quad a,b\in\A
    \end{equation} it is well-defined because of the Cauchy-Schwarz inequality \eqref{eq: CS-in H_Phialpha}. 
    Consider the map 
$J(\alpha):\lambda_\Phi^\alpha(\A)\to \lambda_\Phi^\alpha(\A)$ defined as \begin{equation}
        \label{eq J alpha on lambda A}J(\alpha)\lambda_\Phi^\alpha(a):=\lambda_\Phi^\alpha(\alpha(a)),\quad a\in\A.
    \end{equation}   
    Then $J(\alpha)$  is well-defined in fact if $a\in N_\Phi^\alpha$ then     $\Phi(a,\alpha(a))=\|\lambda_\Phi^\alpha(a)\|_\alpha^2=0$ and  
    \begin{align*}        \|J(\alpha)\lambda_\Phi^\alpha(a)\|_\alpha&=\rip{J(\alpha)\lambda_\Phi^\alpha(a)}{J(\alpha)\lambda_\Phi^\alpha(a)}_\alpha\\&=\Phi(\alpha(a),a)=\Phi(a,\alpha(a))=0.
    \end{align*}

      $J(\alpha)$ fulfills $(j1)-(j3)$  in Definition \ref{defn-J-repre}: $J(\alpha)^2\lambda_\Phi^\alpha(a)=\lambda_\Phi^\alpha(\alpha^2(a))=\lambda_\Phi^\alpha(a)$, for every $a\in\A$,  $$\sip{\lambda_\Phi^\alpha(a)}{J(\alpha)\lambda_\Phi^\alpha(b)}=\Phi(a,\alpha(b))=\Phi(\alpha(a),\alpha^2(b))=\sip{J(\alpha)\lambda_\Phi^\alpha(a)}{\lambda_\Phi^\alpha(b)},
    $$ for all $a,b\in\A$ and  \begin{equation*}
       \rip{\lambda_\Phi^\alpha(a)}{\lambda_\Phi^\alpha(b)}_\alpha=\sip{\lambda_\Phi^\alpha(a)}{J(\alpha)\lambda_\Phi^\alpha(b)},\quad a,b\in\A
    \end{equation*} therefore  $\alpha\in j_{\lambda_\Phi^\alpha(\A)}(\A)$ and, by Remark \ref{rem boundedness of J(alpha)}, $J(\alpha)$ is bounded on $\lambda_\Phi^\alpha(\A)$ with respect to the norm $\|\cdot\|_\alpha$. \\
    
    In what follows 
 $\K_\Phi$ will denote the completion of $\lambda_\Phi^\alpha(\A)$ with respect to $\|\cdot\|_\alpha=\|\cdot\|_{J(\alpha)}$, then $\K_\Phi$ is a Hilbert space with respect to the definite positive  inner product in  \eqref{eq def pos inner prod on lambda A}. 
By Theorem \ref{thm extendability of the inner prod} an indefinite inner product  
    extending $\sip{\cdot}{\cdot}$ can be defined on $\K_\Phi$ (still denoted by  $\sip{\cdot}{\cdot}$, see \eqref{eq indef inner prod on lambda A}) such that $(\K_\Phi, \sip{\cdot}{\cdot})$ is a Krein space.

         \begin{thm} \label{thm: GNS per Phi in Krein} Let $(\A,\A_0)$ be a quasi *-algebra with unit $\id$, $\alpha\in j(\A)$ and
              $\Phi\in\QAal$ 
             be such that $\lambda_\Phi^\alpha(\A_0)$  is a dense subspace of the Krein space $\K_\Phi$. Then, there exist a Krein space $\{\K_\Phi^\alpha,\sip{\cdot}{\cdot}_\pi\}$, a  dense domain $\D_\Phi^\alpha$ of  $\K_\Phi^\alpha$ and a closed  $\alpha$-cyclic $J(\alpha)$-representation $\pi_\Phi^\alpha$ in ${\mathcal L}^\ddagger(\D_\Phi^\alpha,\K_\Phi^\alpha)$ of $(\A,\A_0)$, with cyclic vector $\xi_\Phi^\alpha$  such that \begin{equation*}
\Phi(a,b)=\sip{\pi_\Phi^\alpha(a)\xi_\Phi^\alpha}{\pi_\Phi^\alpha(b)\xi_\Phi^\alpha}_{\pi},\quad \forall a,b\in\A\end{equation*} and  for all $a\in\A$ \begin{equation*}        \pi_\Phi^\alpha(\alpha(a))\eta=J(\alpha)\pi_\Phi^\alpha(a) J(\alpha)\eta,\quad \forall \eta\in\D_\Phi^\alpha.
     \end{equation*}
\end{thm}\begin{proof} Put $\D_{\pi^\circ}=\lambda_\Phi^\alpha(\A_0)$, then $\D_{\pi^\circ}$ is dense in the Hilbert space $\K_{\Phi}^\alpha=\K_\Phi$. 
Put, for any $a\in\A$
$$\pi^\circ(a)\lambda_\Phi^\alpha(x):=\lambda_\Phi^\alpha(ax), \quad x\in \A_0.$$ 
It is well-defined since, if $x\in   N_\Phi^\alpha\cap \A_0$ i.e. $\lambda_\Phi^\alpha(x)=0$ then, by \eqref{eq: property of N},    for any $a\in\A$ $$\Phi(ax,\alpha(y))=\Phi(x,\alpha(a^*y))=0\quad \forall y\in\A_0.$$ Now,
by the density of $\lambda_\Phi^\alpha(\A_0)$  in $\K_\Phi$, for all $b\in\A$ there exists a sequence  $\{y_n\}\subseteq\A_0$, such that 
$\|\lambda_\Phi^\alpha(b)-\lambda_\Phi^\alpha(y_n)\|_\alpha\to 0$ as $n\to\infty$ i.e., $$\Phi(b-y_n,\alpha(b-y_n))\to0$$ as $n\to \infty$. By \eqref{eq for the well-definiteness} it follows  that $\Phi(ax,\alpha(y_n))=\Phi(x,\alpha(a^*y_n))=0$ since $x\in\A_0$. Thus, since $\Phi$ is $\alpha$-positive it is      
  \begin{align*}
    |\Phi(ax,\alpha(b))|&= |\Phi(ax,\alpha(b-y_n))+ \Phi(ax,\alpha(y_n))|\\\nonumber&\leq|\Phi(ax, \alpha(b-y_n))|\\\nonumber&\leq \Phi(ax, \alpha(ax))^{1/2}\Phi(b-y_n, \alpha(b-y_n))^{1/2}\to0,
    \end{align*}\,as $n\to\infty$, for every  $b \in \A$, hence $\Phi(ax,\alpha(b))=0$, for every  $b \in \A$, therefore $ax\in N_\Phi^\alpha$, i.e. $\pi^\circ(a)\lambda_\Phi^\alpha(x)=\lambda_\Phi^\alpha(ax)=0$.

   Consider the  operator $J(\alpha)$ defined in \eqref{eq J alpha on lambda A}, and restrict it to $\D_{\pi^\circ}$. By  $(\alpha.3)$, one has $$J(\alpha):\D_{\pi^\circ}\to\D_{\pi^\circ}$$ hence $\alpha\in j_{\D_{\pi^\circ}}(\A)$. 
   Consider now  the restriction of the indefinite inner
product defined in \eqref{eq indef inner prod on lambda A} to $\D_{\pi}^\circ$. Since
$J(\alpha):\D_{\pi}^\circ\to \D_{\pi}^\circ$ and satisfies
$(j1)$--$(j3)$, Theorem \ref{thm extendability of the inner prod}
implies that the Hilbert space $\K_\Phi$, endowed with the extension of
this indefinite inner product, is a Krein space with fundamental
symmetry $J(\alpha)$.
   By setting  $\xi_\Phi^\alpha=\lambda_\Phi^\alpha(\id)\in\D_{\pi^\circ}$, then $J(\alpha)\xi_\Phi^\alpha=\xi_\Phi^\alpha$ and 
$\xi_\Phi^\alpha$ is such that $\pi^\circ(\A_0)\xi_\Phi^\alpha=\lambda_\Phi^\alpha(\A_0)$  is dense in $\K_\Phi$ in its norm topology $\|\cdot\|_{J(\alpha)}=\|\cdot\|_\alpha$.
Moreover,  for every $a\in\A$, $\eta\in\D_{\pi^\circ}$, $\pi^\circ(\alpha(a))\eta=J(\alpha)\pi^\circ(a)J(\alpha)\eta$, in fact, for every $\lambda_\Phi^\alpha(x)\in \D_{\pi^\circ}$, 

\begin{align*}\pi^\circ(\alpha(a))\lambda_\Phi^\alpha(x)&=\lambda_\Phi^\alpha(\alpha(a)x)=\lambda_\Phi^\alpha(\alpha(a\alpha(x)))\\
&=J(\alpha)\lambda_\Phi^\alpha(a\alpha(x)) =J(\alpha)\pi^\circ(a)\lambda_\Phi^\alpha(\alpha(x))\\&=J(\alpha)\pi^\circ(a)J(\alpha)\lambda_\Phi^\alpha(x).
    \end{align*}
   Before proceeding, observe that, for every $a\in\A$,
\begin{equation}\label{eq hilb adj}
\pi^\circ(a^*)\subseteq \pi^\circ(a)^*.
\end{equation}
 where  $\pi^\circ(a)^*$ simply denotes the operator $\pi^\circ(a)^*_\alpha$ i.e. the Hilbertian adjoint of $\pi^\circ(a)$ with respect to the positive definite inner product $\rip{\cdot}{\cdot}_\alpha$. In fact,  \begin{align*}
    \rip{\pi^\circ(a)\lambda_\Phi^\alpha(x)}{\lambda_\Phi^\alpha(y)}_\alpha&=\Phi(ax,\alpha(y))=\Phi(x,\alpha(a^*y))\\&=\rip{\lambda_\Phi^\alpha(x)}{\pi^\circ(a^*)\lambda_\Phi^\alpha(y)}_\alpha,\qquad\forall x,y\in\A_0.
\end{align*}

Since $\D_\pi^\circ$ is dense in $\K_\Phi^\alpha$, \eqref{eq hilb adj} implies that
$\pi^\circ(a)^*$ is densely defined. Hence $\pi^\circ(a)$ is closable. 
The next step is to prove that the map $\pi^\circ$ is a $J(\alpha)$-representation of $(\A,\A_0)$ in ${\mathcal L}^\ddagger(\D_{\pi^\circ},\K_\pi)$. 
The operators $\pi^\circ(\alpha(a^*))$ and $\pi^\circ(a)^{J(\alpha)}$ coincide on $\D_{\pi^\circ}$. In fact, by Remark \ref{rem relationship bet adjoints}, for every $a\in\A$  it is $\pi^\circ(a)^{J(\alpha)}=J(\alpha)\pi^\circ(a)^*J(\alpha)$ and together with \eqref{eq hilb adj} and the covariance $\pi^\circ(\alpha(a))=J(\alpha)\pi^\circ(a)J(\alpha)$, it yields:
$$\pi^\circ(\alpha(a^*))=J(\alpha)\pi^\circ(a^*)J(\alpha)\subseteq J(\alpha)\pi^\circ(a)^*J(\alpha)=\pi^\circ(a)^{J(\alpha)},$$
therefore the restriction $\pi^\circ(a)^\ddagger$ of $\pi^\circ(a)^{J(\alpha)}$ to $\D_{\pi^\circ}$ coincides with $\pi(\alpha(a^*))$, then \begin{multline*}\pi(a^*)=J(\alpha)\pi(\alpha(a^*))J(\alpha)= J(\alpha)\pi(\alpha(a)^*)J(\alpha)\\\subseteq J(\alpha)\pi(\alpha(a))^*J(\alpha)=\pi(\alpha(a))^{J(\alpha)}\end{multline*} and $\pi(a^*)=\pi(\alpha(a))^{J(\alpha)}{\upharpoonright}_{\D_{\pi^\circ}}$. Hence $(r1)$ is proved.

Observe that for every $a\in\A$ by Lemma \ref{lem 3.1}, since $\pi^\circ(a)$ is closable, one has $$\pi^\circ(a)\subseteq \pi^\circ(a)^{J(\alpha)J(\alpha)}$$ and by $\pi^\circ(a)^\ddagger\subseteq\pi^\circ(a)^{J(\alpha)}$ it is  $$\pi^\circ(a)^{J(\alpha)J(\alpha)}\subseteq (\pi^\circ(a)^\ddagger)^{J(\alpha)},$$ therefore $\pi^\circ(a)\subseteq  (\pi^\circ(a)^\ddagger)^{J(\alpha)}$.
For $a\in\A$ and $ x\in\A_0$,  it is $$\pi^\circ(x)\lambda_\Phi^\alpha(y)=\lambda_\Phi^\alpha(xy)\in\D_{\pi^\circ}\subseteq D((\pi^\circ(a)^\ddagger)^{J(\alpha)}),\quad \forall y\in\A_0.$$ To prove that   
$$\pi^\circ(a)^\ddagger\D_{\pi^\circ}\subset D(\pi^\circ(x)^{J(\alpha)}),$$
notice that, since $\pi^\circ(a)^\ddagger\lambda_\Phi^\alpha(z)=\pi^\circ(\alpha(a^*))\lambda_\Phi^\alpha(z)$ for every $z\in\A_0$, then\begin{align*}
    \sip{\pi^\circ(x)\lambda_\Phi^\alpha(y)}{\pi^\circ(\alpha(a^*))\lambda_\Phi^\alpha(z)}&=\Phi(xy,\alpha(a^*)z)=\Phi(xy,\alpha(a^*\alpha(z)))\\
&=\Phi(\alpha(x)\alpha(y),\alpha(\alpha(a^*)z)))\\
&=\Phi(\alpha(y),\alpha(\alpha(x^*a^*)z)))\\&=\Phi(\alpha(y),\alpha(\alpha(x^*a^*))\alpha(z))\\&=\Phi(\alpha(y),x^*a^*\alpha(z))\\&=\Phi(y,\alpha(x^*a^*)z))\\&=\sip{\lambda_\Phi^\alpha(y)}{\lambda_\Phi^\alpha(\alpha(x^*a^*)z))}\end{align*}hence, there exists $\eta=\lambda_\Phi^\alpha(\alpha(x^*a^*)z))\in\K_\Phi$ such that $$\sip{\pi^\circ(x)\lambda_\Phi^\alpha(y)}{\pi^\circ(\alpha(a^*))\lambda_\Phi^\alpha(z)}=\sip{\lambda_\Phi^\alpha(y)}{\eta}, \quad \forall y\in\A_0$$
therefore $\lambda_\Phi^\alpha(\alpha(x^*a^*)z))\in D(\pi^\circ(x)^{J(\alpha)})$. Then  the two operators $\pi^\circ(a)$ and $\pi^\circ(x)$ can be weakly multiplied. It is 
well-defined because,  
if $y\in  N_\Phi^\alpha\cap \A_0$ (i.e. $\lambda_\Phi^\alpha(y)=0$), then for what it has been proved before $xy\in N_\Phi^\alpha$ and
\begin{align*}
\sip{\pi^\circ(a)\wmult\pi^\circ(x)\lambda_\Phi^\alpha(y)}{\lambda_\Phi^\alpha(z)}&=\sip{(\pi^\circ(a)^\ddag)^{J(\alpha)} \pi^\circ(x) \lambda_\Phi^\alpha(y)}{\lambda_\Phi^\alpha(z)}\\&= \sip{  \lambda_\Phi^\alpha(xy)}{\pi^\circ(a)^\ddag\lambda_\Phi^\alpha(z)}\\&= \sip{  \lambda_\Phi^\alpha(xy)}{\lambda_\Phi^\alpha(\alpha(a^*)z)}=0,\quad\forall z\in\A_0
\end{align*} 
moreover,
\begin{align*}
\pi^\circ(a)\wmult\pi^\circ(x)\lambda_\Phi^\alpha(y)&=(\pi^\circ(a)^\ddagger)^{J(\alpha)}\lambda_\Phi^\alpha(xy)\\&=\pi^\circ(a)\lambda_\Phi^\alpha(xy)=\lambda_\Phi^\alpha(axy)=\pi^\circ(ax)\lambda_\Phi^\alpha(y), \,\,\forall y\in\A_0.
\end{align*}

Now, consider the closed  $J(\alpha)$-representation $\widetilde{\pi}$ defined as $$\widetilde{\pi}(a):=\overline{\pi^\circ(a)}{\upharpoonright}_{\widetilde{\D_{\pi^\circ}}},\quad a\in\A$$ the restriction of the closure $\overline{\pi^\circ(a)}$ of the operator $\pi^\circ(a)$ to the completion $\D_\Phi^\alpha=\widetilde{\D_{\pi^\circ}}$ with respect to the graph topology generated by the family of seminorms \eqref{eq: seminorms} with $J_\alpha=J(\alpha)$. By Proposition \ref{prop closed repr}  $\widetilde{\pi}$ is a closed $J(\alpha)$-representation and  the vector $\xi^\alpha_\Phi$ is $\alpha$-cyclic also for $\widetilde{\pi}$. Thus, after renaming  $\pi_\Phi^\alpha=\widetilde{\pi}$ one obtains $$\Phi(a,b)=\sip{\pi_\Phi^\alpha(a)\xi_\Phi^\alpha}{\pi_\Phi^\alpha(b)\xi_\Phi^\alpha}, \quad \forall a,b\in\A.$$ 
\end{proof}

 \begin{rem}
     Theorem \ref{thm: GNS per Phi in Krein} is genuinely a generalization of the classical GNS-construction for a  {\em positive} sesquilinear form, see Proposition \ref{GNS with a pos l.i. s.f.}, in fact when $\alpha=I$ and $\Phi$ is positive, the construction reduces to
the classical Hilbert-space GNS construction.
 \end{rem}

\section{Abstract conditional expectation for a quasi *-algebra and GNS}\label{sec abstract expect}
In \cite{Sakai} the notion of abstract conditional expectation of a *-algebra $\A_0$ with unit  onto a *-subalgebra with the same unit  has been introduced  and then, in \cite{AO}, it has been exploited to construct a map $\alpha$ on $\A_0$ that plays the role of the automorphism $\alpha$ in the GNS construction for a linear functional on $\A_0$ as in \cite{Jakobczyk86}. In 2009,  Takakura \cite{Tak} introduced the notion of unbounded conditional expectation for partial O$^*$-algebras. In this section the notion of abstract conditional expectation of a quasi *-algebra $(\A,\A_0)$ with unit $\id$ onto a *-invariant subspace $\FF$ containing the same unit 
$\id$ is introduced and then used to construct a linear map $\alpha\in j(\A)$ that, together with an additional condition, allows a GNS construction for suitable sesquilinear forms on $\A$ in a Krein space. \\

\begin{defn}\label{defn ace}
    Let $(\A,\A_0)$ be a quasi *-algebra with unit $\id$.
Let $\FF$ be a *-invariant linear subspace of $\A$ containing $\id$, and 
$
\FF_0:=\FF\cap \A_0
$
be a unital *-subalgebra of $\A_0$ with the same unit. A linear map $P:\A\to \FF$ is called an
{\em abstract conditional expectation of $\A$ onto $\FF$} if:
\begin{itemize}
\item[(i)] $P(\A)=\FF$ and $P(f)=f$, for every $f\in \FF$;
\item[(ii)] $P(a^*)=P(a)^*$, for every $a\in \A$;
\item[(iii)] $P(\A_0)\subseteq \A_0$;
\item[(iv)] $P(xay)=xP(a)y$, for every $a\in \A$ and every
$x,y\in \FF_0$.
\end{itemize}
\end{defn}

As a direct consequence of the Definition \ref{defn ace}, $P^2(a)=P(P(a))=P(a)$ for all $a\in\A$ and $P(\id)=\id$.   Moreover, Definition \ref{defn ace} generalizes the bimodularity of the definition in \cite{AO}, by replacing the subalgebra $\mathfrak{B}$ of the *-algebra $\A_0$,  with $\FF_0 = \FF \cap\A_0$, because, as it is known, in a quasi *-algebra, two arbitrary elements of $\A$ cannot be multiplied. \\

The following is a simple example of abstract conditional expectation of a quasi *-algebra with unit onto a *-invariant linear subspace of $\A$.
\begin{ex}\label{ex ace not admissible}
	Let
	$
	(\A,\A_0)=(L^2([0,1]),L^\infty([0,1]))
	$
	  and 
	$
	\FF=\mathbb{C}\id
	$, with $\id=1$ the constant function on $[0,1]$. 
	Then,
	$$
	\FF_0=\FF\cap \A_0=\mathbb{C}\id,
	$$
	which is a *-subalgebra of $\A_0$. Define $P:\A\to\FF$ as
	$$
	P(f):=\left(\int_0^1 f(t)\,dt\right)\id,\qquad f\in L^2([0,1]).
	$$
	Then $P:\A\to \FF$ is an abstract conditional expectation of $\A$ onto $\FF$.
		\end{ex}
\medskip

In analogy with \cite{AO}, the following class of abstract conditional expectations is considered.
\begin{defn}  Let $(\A,\A_0)$ be a quasi *-algebra with unit $\id$.
Let $\FF$ be a *-invariant linear subspace of $\A$ containing $\id$.
   An  abstract conditional expectation $P$ of $\A$ onto $\FF$ is called {\em admissible} if the map
\begin{equation}
    \label{eq: alpha_AS}
\alpha_P:=2P-I
\end{equation}
belongs to $j(\A)$, with $I$ the identity map on $\A$.
\end{defn}

\begin{rem}\label{rem fixed points}
Let $P:\A\to \FF$ be an abstract conditional expectation of $\A$ onto $\FF$,
and put
$$
\alpha_P:=2P-I.
$$
Then
$$
F_{\alpha_P}(\A)=\FF.
$$
Indeed, for $a\in \A$,
$$
\alpha_P(a)=a
$$
if and only if
$$
2P(a)-a=a,
$$
which is equivalent to
$$
P(a)=a.
$$
Since $P$ is the identity on $\FF$ and $P(\A)=\FF$, this holds if and only if
$a\in \FF$. In particular, if $P$ is admissible, then
$$
F_{\alpha_P}(\A_0)
=
F_{\alpha_P}(\A)\cap \A_0
=
\FF\cap \A_0
=
\FF_0.
$$
\end{rem}
This further condition on $P$  is an additional
compatibility one between the abstract conditional expectation
and the quasi *-algebra structure, as will be shown below.

	\begin{rem}\label{rem P non admissi} Admissibility is a genuine additional assumption for $P$,
needed in order to associate with $P$ a map $\alpha_P=2P-I$ belonging
to $j(\A)$. In fact, consider the quasi *-algebra and the abstract conditional expectation of the Example \ref{ex ace not admissible} and the linear map 
		$
	\alpha_P:=2P-I.
	$
	Then
	$$
	\alpha_P(f)=2\left(\int_0^1 f(t)\,dt\right)\id-f.
	$$

	It is $
	\alpha_P\notin j(\A)
	$. Indeed, take $
	f(t)=t$ and $h(t)=t$, then 	$(fh)(t)=t^2$,
	hence
	$$ \alpha_P(fh)=2\int_0^1 t^2\,dt - t^2	=\frac{2}{3}-t^2 $$
	but on the other hand,
	$ \alpha_P(f)	=	1-t$ and $	\alpha_P(h)	=
	1-t	$,
	i.e.,
	$$ 	\alpha_P(f)\alpha_P(h)=(1-t)^2. 	$$
	Therefore,
	$$
	\alpha_P(fh)=\frac{2}{3}-t^2\neq (1-t)^2= \alpha_P(f)\alpha_P(h).
	$$
	\end{rem}

\bigskip
However, there are examples of admissible abstract conditional expectation. One of the simplest of them is shown in Example \ref{ex adm ace1} below.

\bigskip

Let $P:\A\to \FF$ be an abstract conditional expectation of $\A$ onto
$\FF$, and put
$$
\MM:=\ker P.
$$
Then
$$
\A=\FF\oplus \MM.
$$
Indeed, since $P$ is the identity on $\FF$, every $a\in\A$ can be written
uniquely as
$$
a=P(a)+(a-P(a)),
$$
with $P(a)\in\FF$ and $a-P(a)\in\MM$. Moreover, since
$P(\A_0)\subseteq \A_0$, one has
$$
\A_0=\FF_0\oplus \MM_0,
\qquad
\FF_0:=\FF\cap \A_0,\quad \MM_0:=\MM\cap \A_0.
$$

With respect to the map
$
\alpha_P:=2P-I
$,
the space $\FF$ is the $+1$-eigenspace and $\MM$ is the $-1$-eigenspace.
Thus $\FF_0$ and $\MM_0$ may be regarded, respectively, as the even and
odd parts of $\A_0$ with respect to $\alpha_P$, even when
$\alpha_P$ does not belong to $j(\A)$.\\

    The right multiplication by
elements of $\FF_0$ always preserves the decomposition
$$
\A=\FF\oplus\MM.
$$
Indeed, if $u\in\FF_0$ and $a=f+m\in\FF\oplus\MM$, then
$$
P(au)=P(a)u=fu.
$$
Hence
$$
au-fu=mu\in\MM.
$$
Thus
$$
\FF\FF_0\subseteq\FF,
\qquad
\MM\FF_0\subseteq\MM.
$$
Therefore, in order to characterize admissibility, it remains only to
control right multiplication by elements of $\MM_0$.   
    In fact, the following characterization can be given.
\begin{prop}\label{prop charact of admiss P}
Let $P:\A\to \FF$ be an abstract conditional expectation of $\A$ onto
$\FF$, and put
$$
\MM:=\ker P,\qquad
\FF_0:=\FF\cap \A_0,\qquad
\MM_0:=\MM\cap \A_0.
$$
Then the following statements are equivalent:
\begin{itemize}
\item[(i)] $P$ is admissible;
\item[(ii)]
$
\FF\MM_0\subseteq \MM
\quad\text{and}\quad
\MM\MM_0\subseteq \FF,
$
where products are understood as right module products.
\end{itemize}
\end{prop}
\begin{proof}
First observe that, for every abstract conditional expectation $P$, the
map
$
\alpha_P:=2P-I
$ satisfies properties $(\alpha.1)$, $(\alpha.2)$, $(\alpha.4)$ and moreover $
\alpha_P(\A_0)\subseteq \A_0$.
Hence $P$ is admissible if and only if
$$
\alpha_P(ax)=\alpha_P(a)\alpha_P(x),
\qquad a\in\A,\ x\in\A_0.
$$

Let $
a=f+m\in \FF\oplus\MM$, $x=u+v\in \FF_0\oplus\MM_0$, then $\alpha_P(a)=f-m$,
and $\alpha_P(x)=u-v$.

Let $v\in\MM_0$. The condition
$$
\alpha_P(av)=\alpha_P(a)\alpha_P(v)
$$
becomes
$$
\alpha_P((f+m)v)=(f-m)(-v),
$$
that is,
$$
\alpha_P(fv+mv)=-fv+mv.
$$
Since $\alpha_P$ is the identity on $\FF$ and minus the identity on
$\MM$, this equality holds for every $f\in\FF$, $m\in\MM$ and
$v\in\MM_0$ if and only if
$$
fv\in\MM
\quad\text{and}\quad
mv\in\FF.
$$
Therefore $P$ is admissible if and only if
$
\FF\MM_0\subseteq\MM$ and $
\MM\MM_0\subseteq\FF
$.
\end{proof}

\begin{ex}\label{ex adm ace1}
   Let
	$(\A,\A_0)=(L^2([0,1]),L^\infty([0,1])).$
	Again, $(\A,\A_0)$ is a quasi *-algebra with pointwise multiplication
	whenever one of the two factors belongs to $\A_0$, involution
	$f^*=\overline f$, and unit $\id=1$.
	Let $R$ be the reflection operator
	$$
	(Rf)(t)=f(1-t),
	\qquad f\in L^2([0,1]).
	$$
	Define
	$$
	P=\frac{I+R}{2},
	$$
	that is,
	$$
	P(f)(t)=\frac{f(t)+f(1-t)}{2}.
	$$
	The range of $P$ is
	$$
	\FF=\{f\in L^2([0,1]): f(t)=f(1-t)\ \text{a.e.}\}.
	$$In fact, $P(f)(1-t)=\frac{f(t)+f(1-t)}{2}=P(f)(t)$, hence $P(\A)\subseteq\FF$. On the other hand, if $g\in\FF$, then $g(t)=g(1-t)$ a.e., therefore $$P(g)(t)=\frac{g(t)+g(1-t)}{2}=g(t),$$ hence $\FF\subseteq P(\A)$. 
	Moreover,
	$$
	\FF_0=\FF\cap \A_0
	=
	\{h\in L^\infty([0,1]): h(t)=h(1-t)\ \text{a.e.}\}, 
	$$ 
	which is a *-subalgebra of $\A_0$. It is easily shown that $P$ is admissible. It is $L^2([0,1])=\FF\oplus\MM$ with
$$
\MM:=\ker P
=
\{f\in L^2([0,1]): f(1-t)=-f(t)\ \text{a.e.}\},
$$

$$
\FF
=
\{f\in L^2([0,1]): f(1-t)=f(t)\ \text{a.e.}\},
$$
and
$$
\FF_0=\FF\cap L^\infty([0,1]),
\qquad
\MM_0=\MM\cap L^\infty([0,1]).
$$

Let $v\in\MM_0$. If $f\in\FF$, then
$$
(fv)(1-t)
=
f(1-t)v(1-t)
=
f(t)(-v(t))
=
-(fv)(t),
$$
hence $fv\in\MM$. Therefore $
\FF\MM_0\subseteq\MM.
$
On the other hand, if $m\in\MM$, then
$$
(mv)(1-t)
=
m(1-t)v(1-t)
=
(-m(t))(-v(t))
=
m(t)v(t),
$$
hence $mv\in\FF$. Therefore $
\MM\MM_0\subseteq\FF
$. By Proposition \ref{prop charact of admiss P}, $P$ is admissible.
\end{ex}

\begin{defn}\label{defn: Phi}
   Let $(\A,\A_0)$ be a quasi *-algebra with unit $\id$ and let  $P:\A\to\FF$ be an
 abstract conditional expectation of $\A$ onto $\FF$. A hermitian (i.e. $\Phi(b,a)=\overline{\Phi(a,b)}$, for all $a,b\in\A$) sesquilinear form $\Phi$ on $\A$ is called a
$P$-form on $\A$ if \begin{itemize}
   \item[(P1)]  $\Phi(P(a),b)=\Phi(a,P(b))=\Phi(P(a),P(b)), \,\, \forall  a,b \in \A$,
    \item[(P2)] $\Phi(a,a)\leq 2 \Phi(P(a),P(a))$, for all $a\in\A$,
    \item[(P3)]  
for every $a\in \A$ and $x,y\in \A_0$,
$$
2\Phi(ax,P(y))-\Phi(ax,y)
=
2\Phi(x,P(a^*y))-\Phi(x,a^*y).
$$
    \end{itemize}
\end{defn} 

\begin{rem}
Since $\Phi$ is hermitian, both $\Phi(a,a)$ and
$\Phi(P(a),P(a))$ belong to $\mathbb{R}$ for every $a\in \A$.
In general, however, they do not have a fixed sign.

Moreover, by $(P1)$, for every $a\in \A$ and every $x\in \A_0$ one has
$$\Phi(x,P(a))=\Phi(P(x),P(a))=\Phi(P(x),a). $$
In particular, if $x\in \FF_0=\FF\cap \A_0$, then $P(x)=x$, and therefore
$$ \Phi(x,P(a))=\Phi(x,a), \qquad a\in \A,\ x\in \FF_0. $$
\end{rem}

\begin{ex}
	The abstract conditional expectation $P$ of $\A$ 	onto $\FF$ defined in Example \ref{ex ace not admissible} is not admissible, see Remark \ref{rem P non admissi}. 
	Nevertheless, one can define a $P$-form associated with this non
	admissible expectation; consider, e.g.,
		$$
	\Phi(f,g)
	=
	2m(f)\overline{m(g)}
	-
	\int_0^1 f(t)\overline{g(t)}\,dt,
	$$ with $$
	m(f):=\int_0^1 f(t)\,dt,
	\qquad f\in L^2([0,1]).
	$$
	
It is worth noting that $\Phi$ is an {\em indefinite} $P$-form.
\end{ex}

On the contrary, there are examples of $P$-forms such that  the
	corresponding abstract conditional expectation $P$ is admissible.
\begin{ex}\label{ex P-form}

	Let $\Phi\in {\mathcal Q}_{L^\infty([0,1])}^\alpha(L^2([0,1]))$

	$$
	\Phi(f,g)=\int_0^1 f(t)\overline{g(1-t)}\,dt
	$$ 
    defined in Example \ref{ex form in Qaa}, with $\alpha(f)(t)=Rf(t)=f(1-t)$, $f\in L^2([0,1])$, as in Example \ref{ex alpha}. Then  $\Phi$
   is also a $P$-form with $P=\displaystyle \frac{I+R}{2}$, as will follow from Proposition \ref{prop charact admiss Palpha} below with $\alpha=R$.
	\end{ex}

\begin{lem}\label{lem P then alpha left-inv alpha pos sesq form}
     Let $(\A,\A_0)$ be a quasi *-algebra with unit $\id$, $\FF$ be a *-invariant subspace of $\A$ containing $\id$, $\FF_0=\FF\cap\A_0$ be a *-subalgebra of $\A_0$ and  $P:\A\to\FF$ an admissible
 abstract conditional expectation of $\A$ onto $\FF$. If $\Phi$ is a P-form  on $\A$, then  $\Phi\in {\mathcal Q}_{\Ao}^{\alpha_P}(\A)$ with $\alpha_P$ defined as in \eqref{eq: alpha_AS}.
\end{lem}
\begin{proof}
    $\Phi$ is  ${\alpha_P}$-positive, in fact, \begin{align*}
        \Phi(a, {\alpha_P}(a))&=\Phi(a, 2P(a)-a)\\&=2\Phi(a, P(a))-\Phi(a,a)\\&=2\Phi(P(a), P(a))-\Phi(a,a)\geq0,\quad \forall a\in\A;
    \end{align*}  besides,   $\Phi$ is clearly ${\alpha_P}$-left-invariant by  $(P3)$, indeed, since $\alpha_P=2P-I$, condition $(P3)$ gives
$$
\Phi(ax,\alpha_P(y))
=
\Phi(x,\alpha_P(a^*y)),
\qquad a\in\A,\ x,y\in\A_0;
$$   
finally, \begin{align*}\Phi({\alpha_P}(a),{\alpha_P}(b))&=\Phi(2P(a)-a,2P(b)-b)\\&=4\Phi(P(a),P(b))\\&\quad-2\Phi(P(a),b)-2\Phi(a,P(b))+\Phi(a,b)\\&=\Phi(a,b), \quad \forall a,b\in\A.\end{align*} 
\end{proof}
  
           The previous discussion associates with every admissible abstract
conditional expectation $P$ an element $\alpha_P=2P-I$ of $j(\A)$.
Conversely, every $\alpha\in j(\A)$ gives rise to a canonical admissible
abstract conditional expectation, namely
\[
P_\alpha=\frac{I+\alpha}{2}.
\]
This correspondence also allows us to reformulate the class
$\QAal$ in terms of $P_\alpha$-forms.
           \begin{prop}\label{prop charact admiss Palpha}Let $(\A,\A_0)$ be a quasi *-algebra with unit $\id$, let
$\alpha:\A\to\A$ be a linear map.
Then the following statements are equivalent:
\begin{itemize}
\item[(i)] $\alpha\in j(\A)$;\\
\item[(ii)] $P_\alpha:=\displaystyle \frac{I+\alpha}{2}$ is an admissible abstract conditional expectation
of $\A$ onto $F_\alpha(\A)$ and
$
\alpha_{P_\alpha}=2P_\alpha-I=\alpha.
$
\end{itemize}
Moreover, if these equivalent conditions hold, then, for every
sesquilinear form $\Phi$ on $\A$,
$$\Phi\in \QAal
\quad\mbox{ if and only if }\quad \Phi \mbox{ is a $P_\alpha$-form}.$$

\end{prop}

\begin{proof}
Assume first that $\alpha\in j(\A)$. By Remark \ref{rem fixed point of alpha in JA}, $F_\alpha(\A)$ is a *-invariant linear subspace
of $\A$ containing $\id$, and
$$
F_\alpha(\A_0)=F_\alpha(\A)\cap \A_0
$$
is a unital *-subalgebra of $\A_0$. Since $\alpha^2=I$, the map
$P_\alpha$ is a projection. Moreover,
$$
P_\alpha(a)=a
$$
if and only if
$$
\alpha(a)=a.
$$
Hence
$$
P_\alpha(\A)=F_\alpha(\A),
$$
and $P_\alpha$ is the identity on $F_\alpha(\A)$. Since $\alpha(\id)=\id$,
one has $P_\alpha(\id)=\id$. Since $\alpha$ commutes with the involution,
$$
P_\alpha(a^*)=P_\alpha(a)^*,
\qquad a\in\A.
$$
Moreover, since $\alpha(\A_0)\subseteq \A_0$, one has
$$
P_\alpha(\A_0)\subseteq \A_0.
$$
Finally, if $x,y\in F_\alpha(\A_0)$ and $a\in\A$, then
$$
\alpha(xay)=x\alpha(a)y.
$$
Therefore
$$
P_\alpha(xay)
=
\frac{xay+x\alpha(a)y}{2}
=
xP_\alpha(a)y.
$$
Thus $P_\alpha$ is an abstract conditional expectation of $\A$ onto
$F_\alpha(\A)$. Since
$$
\alpha_{P_\alpha}=2P_\alpha-I=\alpha,
$$
and $\alpha\in j(\A)$, $P_\alpha$ is admissible.

Conversely, suppose that $P_\alpha$ is admissible and that
$\alpha_{P_\alpha}=\alpha$; by definition
$
\alpha_{P_\alpha}
\alpha\in j(\A).
$
This proves the equivalence.\\

Assume now that these equivalent conditions hold. Since
$$
\alpha=2P_\alpha-I,
$$
condition $(P3)$ is exactly the $\alpha$-left-invariance condition
$(q.2)$.

If $\Phi\in \QAal$, then $\Phi$ is hermitian, because the
positive form
$$(a,b)\mapsto \Phi(a,\alpha(b)),\quad a,b\in\A$$
is hermitian and the covariance gives the hermitianity of $\Phi$.
Let
$$a=a_++a_-,\qquad a_+=P_\alpha(a),\quad a_-=(I-P_\alpha)(a),$$
and similarly
$$ b=b_++b_-,\qquad b_+=P_\alpha(b),\quad b_-=(I-P_\alpha)(b).$$
Then
$$
\alpha(a_+)=a_+,\qquad \alpha(a_-)=-a_-,
$$
and similarly for $b_+$ and $b_-$. By the covariance of $\Phi$, the
even and odd parts are $\Phi$-orthogonal. Indeed,
$$\Phi(a_+,b_-)
= \Phi(\alpha(a_+),\alpha(b_-))=\Phi(a_+,-b_-)=-\Phi(a_+,b_-),$$
hence $\Phi(a_+,b_-)=0$. Similarly,
$\Phi(a_-,b_+)=0.$\\
Therefore,
$$\Phi(P_\alpha(a),b)=\Phi(a_+,b_++b_-)=\Phi(a_+,b_+)=
\Phi(P_\alpha(a),P_\alpha(b)),$$
and
$$\Phi(a,P_\alpha(b))=\Phi(a_++a_-,b_+)=\Phi(a_+,b_+)
=\Phi(P_\alpha(a),P(b)).$$
Thus $(P1)$ follows.

Moreover, by $\alpha$-positivity,
$$
0\leq \Phi(a,\alpha(a))
=
\Phi(a_++a_-,a_+-a_-)
=
\Phi(a_+,a_+)-\Phi(a_-,a_-).
$$
Thus
$$
\Phi(a_-,a_-)\leq \Phi(a_+,a_+),
$$
and therefore
$$
\Phi(a,a)
=
\Phi(a_+,a_+)+\Phi(a_-,a_-)
\leq
2\Phi(a_+,a_+)
=
2\Phi(P_\alpha(a),P_\alpha(a)).
$$
So $(P2)$ holds. Hence $\Phi$ is a $P_\alpha$-form.

Conversely, if $\Phi$ is a $P_\alpha$-form, then $
\Phi\in \mathcal{Q}^{\alpha_{P_\alpha}}_{\A_0}(\A)
$, by Lemma
\ref{lem P then alpha left-inv alpha pos sesq form}.
Since $\alpha_{P_\alpha}=\alpha$, then
$
\Phi\in \QAal.
$
\end{proof}
\medskip

   \begin{thm} \label{thm: GNS for P form Phi in Krein} Let $(\A,\A_0)$ be a quasi *-algebra with unit $\id$ and   $P:\A\to\FF$ be an admissible
 abstract conditional expectation of $\A$ onto $\FF$ and let $\Phi$ be a $P$-form on $\A$  such that $\lambda_\Phi^{\alpha_P}(\A_0)$ is a dense subspace of $\K_\Phi^{\alpha_P}$. Then, there exists an indefinite inner product $\sip{\cdot}{\cdot}_\pi$ such that $\{\K_\Phi^{\alpha_P},\sip{\cdot}{\cdot}_\pi\}$ is a Krein space, a  dense domain $\D_\Phi^{\alpha_P}$ of the Hilbert space $\K_\Phi^{\alpha_P}$ and a closed $\alpha_P$-cyclic $J({\alpha_P})$-representation $\pi_\Phi^{\alpha_P}$ in ${\mathcal L}^\ddagger(\D_\Phi^{\alpha_P},\K_\Phi^{\alpha_P})$ of $(\A,\A_0)$, with cyclic vector $\xi_\Phi^{\alpha_P}$  such that \begin{equation*}
\Phi(a,b)=\sip{\pi_\Phi^{\alpha_P}(a)\xi_\Phi^{\alpha_P}}{\pi_\Phi^{\alpha_P}(b)\xi_\Phi^{\alpha_P}}_{\pi},\quad \forall a,b\in\A\end{equation*}and  for all $a\in\A$ \begin{equation*}        \pi_\Phi^{\alpha_P}({\alpha_P}(a))\eta=J({\alpha_P})\pi_\Phi^{\alpha_P}(a) J({\alpha_P})\eta,\quad \forall \eta\in\D_\Phi^{\alpha_P}.
     \end{equation*}  \end{thm} 
\begin{proof} Since $P$ is admissible, then $\alpha_P\in j(\A)$.
    The thesis follows from  Lemma \ref{lem P then alpha left-inv alpha pos sesq form} and Theorem \ref{thm: GNS per Phi in Krein} applied to $\alpha=\alpha_P$. 
\end{proof}

\begin{ex}
Let $\Phi$ be the $P$-form in Example \ref{ex P-form}. Then,
$$
\Phi(f,\alpha_P(f))
=
\Phi(f,Rf)
=
\int_0^1 |f(t)|^2\,dt.
$$
Hence the positive Hilbert norm associated with $\Phi$ and $\alpha_P$ (i.e. the one associated to the inner product in \eqref{eq def pos inner prod on lambda A})
is exactly the usual $L^2$-norm. Therefore
$$
\K^{\alpha_P}_\Phi=L^2([0,1]).
$$
Since
$$
\A_0=L^\infty([0,1])
$$
is dense in $L^2([0,1])$, one has
$$
\lambda^{\alpha_P}_\Phi(\A_0)
\quad\text{dense in}\quad
\K^{\alpha_P}_\Phi.
$$

Consequently, by Theorem \ref{thm: GNS for P form Phi in Krein}, the GNS-like construction yields a closed
$\alpha_P$-cyclic $J(\alpha_P)$-representation in a Krein space.\end{ex}

\section*{Concluding remarks}

Sesquilinear forms on quasi *-algebras
satisfying a suitable $\alpha$-\,\!\! positivity condition have been proven to give rise to closed
$\alpha$-cyclic $J(\alpha)$-\,\!\! representations in Krein spaces (under an additional condition of density). The
construction does not require the sesquilinear form $\Phi$ itself to be
indefinite. What is needed is the positivity of the associated form
$$
(a,b)\mapsto \Phi(a,\alpha(b)),
\qquad a,b\in\A.
$$
Thus the framework includes genuinely indefinite forms, which are the
main motivation for the present construction, but also covers the case
where the ``indefinite'' inner product obtained from $\Phi$ is in fact
positive definite. In the latter situation, the Krein-space structure
reduces to the Hilbert-space case.

Section \ref{sec abstract expect} shows that abstract conditional expectations provide a
concrete mechanism for constructing linear  maps in $j(\A)$
and, consequently, sesquilinear forms to which the GNS-like
construction applies. More precisely, if $P$ is an admissible abstract
conditional expectation, then
$$
\alpha_P=2P-I
$$
belongs to $j(\A)$, and every $P$-form belongs to
$\mathcal{Q}^{\alpha_P}_{\A_0}(\A)$.

Conversely, every $\alpha\in j(\A)$ gives rise to the canonical
admissible abstract conditional expectation
$$
P_\alpha=\frac{I+\alpha}{2}
$$
onto the fixed-point space $F_\alpha(\A)$. In this case the two
languages are equivalent: a sesquilinear form $\Phi$ belongs to
$\QAal$ if and only if it is a $P_\alpha$-form. Thus
$P$-forms may be viewed as a natural conditional-expectation
reformulation of the $\alpha$-positive forms used in the GNS-like
construction.

This also shows that the distinction between abstract conditional
expectations and admissible ones is essential in the quasi *-algebra
setting: admissibility is precisely the condition that allows the map
$\alpha_P=2P-I$ to interact correctly with the quasi *-algebra
structure.\bigskip\\

\noindent{\bf{Acknowledgements:}} The author is thankful to Prof. C. Trapani for his insightful comments to this work and  acknowledges that it has been carried out within the activities of Gruppo UMI Teoria della Approssimazione e Applicazioni and of GNAMPA of the INdAM.

\bibliographystyle{amsplain}

 \end{document}